\documentclass[11pt]{amsart}
\usepackage{geometry,amsrefs}                
\usepackage{changes, lineno}
\usepackage{graphicx,color}
\usepackage{amssymb}
\usepackage{epstopdf}
\newtheorem{prop}{Proposition}

\newtheorem{theo}{Theorem}
\newcommand{\I}{\operatorname{i}}
\newcommand{\C}{\mathbb C}

\newcommand{\e}{\operatorname{e}}
\renewcommand{\Re}{\operatorname{Re}}
\renewcommand{\Im}{\operatorname{Im}}

\title{On the classificarion of 3-dimensional spherical Sasakian manifolds. }
\author{Daniel Sykes}
\address{D.S. University of New England, School of Science and Technology, Armidale NSW 2351, Australia}
\email{dsykes4@myune.edu.au}
\author{Gerd Schmalz}
\address{G.S. University of New England, School of Science and Technology, Armidale NSW 2351, Australia}
\email{schmalz@une.edu.au}
\author{Vladimir Ezhov}
\address{V.E. Flinders University, College of Science and Engineering, 1284 South Road, Tonsley  SA 5042, Australia; 
MSU, Faculty of Mechanics and Mathematics, Leninskiye Gory 1,  Moscow GSP-1 119991,  Russia}
\email{vladimir.ejov@flinders.edu.au}

\begin{document}
\definechangesauthor[name={GS}, color=red]{GS}
\definechangesauthor[name={D}, color=blue]{D}
\definechangesauthor[name={VE}, color=orange]{VE}
\maketitle

\begin{abstract}
   In this article we consider spherical hypersurfaces in $\C^2$  
    with a fixed Reeb vector field as 3-dimensional Sasakian manifolds. We establish the correspondence between three different sets of parameters, namely, those arising from representing the Reeb vector field as an automorphism of the Heisenberg sphere, the parameters used in Stanton’s description of rigid spheres, and the parameters arising from the rigid normal forms. We also geometrically describe the moduli space for rigid spheres, and provide geometric distinction between Stanton’s hypersurfaces and those found in \cite{MR3515414}.  Finally, we determine the  Sasakian automorphism groups of the rigid spheres and detect the homogeneous Sasakian manifolds among them.
\end{abstract}

\section{Introduction}
A spherical 3-dimensional CR manifold can be locally represented as the Heisenberg sphere 
$$S\colon v=|z|^2,$$
where $(z,w=u+\I v)$ are coordinates in $\mathbb C^2$. Such spherical CR manifold can be turned into a Sasakian manifold by choosing an infinitesimal automorphism $Z$ that is transversal to the complex tangent spaces at each point and making $Z$ the Reeb vector field. Such infinitesimal automorphisms form an 8-parametric family (cf. \eqref{ia}) and different choices of $Z$ may or may not result in an equivalent Sasakian structure. The problem of the classification of the (local) Sasakian structures obtained in this way is closely related to the classification of rigid spheres initiated by N. Stanton in  \cite{Sta90} and completed by the second and third author in \cite{MR3515414}. Recall that a rigid sphere is a spherical CR manifold of the form
$$v=\psi(z),$$  
with $\psi(0)=d\psi(0)=0$. Here $\psi(z)$ can be just $|z|^2$, as above, but it can also be more complicated, e.g. $\arcsin |z|^2$, or any of the implicit functions determined by formula  \eqref{sef}. The absence of $u$ in the equation is equivalent to the vector field $Z_0=\frac{\partial}{\partial u}=\Re \frac{\partial}{\partial w}$ being an infinitesimal automorphism and thus endows the manifold with a Sasakian structure, so that $\frac{\partial}{\partial u}$ becomes the Reeb vector field. The Sasakian manifold $S$ with different choices of Reeb vector fields $Z$ can now be treated as different choices of rigid spherical CR manifolds $M$ with Reeb vector field $Z_0$, since any infinitesimal automorphism $Z$  can be mapped to  $Z_0$ by a local holomorphic coordinate change. Therefore, the classification of local Sasakian  structures on $S$ reduces to the classification of rigid spheres, and more precisely, to the classification of rigid spheres in Stanton's rigid normal form.

Although this classification problem has been solved, the exact relation between different sets of parameters involved in that solution is not entirely clear. This problem was raised by Isaev and Merker in \cite{MR4021084}. More precisely, there are three sets of parameters  
\begin{align}\label{set1}
\tau, \rho \in \mathbb R, \quad a\in\mathbb C\\\label{set2}
\theta, r\in \mathbb R, \quad b\in \mathbb C\\\label{set3}
\theta, r, \phi \in \mathbb R, \quad a\in \mathbb C.
\end{align}
The parameters $(\tau,a,\rho)$ from set \eqref{set1} correspond to the choice of the infinitesimal automorphisms
$$Z= (\I \tau z+ aw +2\I \bar{a} z^2+ \rho zw)\frac{\partial}{\partial z}+ (1+2\I \bar{a}zw + \rho w^2)\frac{\partial}{\partial w},$$
which, in fact, lead to the inequivalent Sasakian structures on $S$ modulo a $\mathbb C^*$ action.

It was shown in \cite{MR3515414} that there is a one-to-one correspondence between the set \eqref{set1} and the power series of rigid spheres in normal form. The parameters from set \eqref{set2} were used by Stanton and correspond to an incomplete list of closed-form expressions for rigid spheres. 
The parameters from set \eqref{set3} have been used in  \cite{MR3515414} to derive a complete list of closed-form expressions for rigid spheres. The auxiliary real parameter $\phi$ from set \eqref{set3} satisfies the algebraic equation
\begin{equation}\label{algseta}
|a|^2= \phi  \left((\theta -2 \phi )^2+r^2\right).
\end{equation}
Below we display the closed form expression of a general rigid sphere from \cite{MR3515414} \footnote{We have corrected here two minor typos in the formula obtained in \cite{MR3515414}.}:
\begin{equation}\label{sef}
(1-4\phi |z|^2) \frac{\sin 2rv}{2r} - \e^{-2\theta v}|z|^2-(\phi - \bar{a}z -a \bar{z}+4\phi(\phi-\theta)|z|^2) \frac{\e^{-2\theta v}-\cos 2rv + \frac{\theta\sin 2rv}{r}}{r^2+\theta^2} =0.
\end{equation}
Theorem \ref{th1} below clarifies the exact relation between the parameters and the equivalence classes of rigid spheres.
\begin{theo}\label{th1}
For each pair $(\tau,\rho)\in\mathbb R^2$ such that $\tau^4+\rho^2\le 1$ there exists exactly one class of equivalent rigid spheres, which is represented by equation \eqref{sef} with parameters 
\begin{align}\label{paratheta}
\theta=& \tau + 3\phi\\\nonumber
r^2=& -\rho + (2\tau+3\phi)\phi\\\nonumber
a=&(1-\tau^4 - \rho^2)^{3/8}
\end{align}
and any real solution  $\phi$ of
\begin{equation}\label{algset}
|a|^2= 4\phi^3+ 4\tau\phi^2+ (\tau^2-\rho)\phi.
\end{equation}
Together with the Heisenberg sphere this covers all classes of equivalent rigid spheres.
\end{theo}

Since the equation \eqref{algseta} can have 3 real solutions the choice of $\phi$ may be ambiguous. In Proposition \ref{branch} we answer a question asked by Isaev, by showing that there is no continuous branch of real solutions and that, coincidentally, ambiguity arises exactly in the cases that did not occur in Stanton's list, i.e. the cases where the parameters  \eqref{set2} do not exist. We also find a parametric representation of the real discriminant curve of the cubic \eqref{algset}.

Finally, we compute the symmetry algebra for each equivalence class of spherical Sasakian manifolds. We show

\begin{theo}\label{th2}
Let $M$ be a 3-dimensional Sasakian manifold with underlying spherical CR manifold. Let (without loss of generality)
$$Z=(\I \tau z+ aw +2\I \bar{a} z^2 + \rho zw )\frac{\partial}{\partial z} + (1+2\I \bar{a} zw + \rho w^2) \frac{\partial}{\partial w}$$
be the Reeb vector field. Then the Lie algebra of infinitesimal Sasaki automorphisms is 
\begin{enumerate}
\item 4-dimensional in the case when $\tau=0, a=0, \rho=0$. This is the Heisenberg sphere $v=|z|^2$ with the affine automorphisms 
$$(z,u)\mapsto (\e^{\I\phi}(z+p), u+q-2\Im z\bar{p}),$$
where $p\in\mathbb C$, $q\in\mathbb R$, $\phi\in[0,2\pi)$.
\item 4-dimensional in the case  $a=0$, $\rho=\tau^2$ and $\tau>0$. This corresponds to the rigid sphere  $v=\log(1+|z|^2)$ and can be globally realised as the round sphere $|z|^2+|w|^2=1$ in $\mathbb C^2$ with its natural Sasakian structure.  The Sasakian automorphisms are induced by the unitary transformations in $\mathbb C^2$.
\item 4-dimensional in the case  $a=0$, $\rho=\tau^2$ and $\tau<0$. This corresponds to the rigid sphere  $v=-\log(1-|z|^2)$ and can be globally realised as the hyperboloid $|z|^2-|w|^2=1$ in $\mathbb C^2$ with its natural Sasakian structure.  The Sasakian automorphisms are induced by the pseudounitary transformations in $\mathbb C^2$.
\item 2-dimensional in all other cases.
\end{enumerate} 
\end{theo}

\section{Rigid normal form and equivalence}
For the convenience of the reader we recall here the notions of rigid normal form and equivalence.

A real hypersurface $M \subset \mathbb C^2$ is called rigid if it is invariant with respect to parallel displacement along a real line that is transversal to the complex tangent plane (see \cite{BER}). Geometrically this means that $M$ is a cylinder over a two-dimensional surface in a three-dimensional space transversal to the line of displacement. We choose coordinates $z,w=u+\I v$ in $\mathbb C^2$, where the line of displacement is the $u$-axis and the hypersurface $M$ has the (local) graph equation 
$$v=\psi(z)$$  
with $\psi(0)=d\psi(0)=0$. From now on we assume that $M$ is real-analytic and Levi-nondegenerate, i.e. $\frac{\partial^2\psi}{\partial z\partial \bar{z}}\neq 0$.

N. Stanton \cite{Sta90} has shown that one can choose local holomorphic coordinates in $\mathbb C^2$ such that the defining function of $M$ takes the following rigid normal form
\begin{equation}\label{rnf}
v=|z|^2 + \sum_{k,\ell\ge 2}\gamma_{k\ell}z^k\bar{z}^\ell,
\end{equation}
where $\gamma_{k\ell}$ with $k,\ell \ge 2$ are constants that satisfy $\gamma_{\ell k}=\overline{\gamma_{k\ell}}$.

Stanton noticed that the Heisenberg sphere $v=|z|^2$ can have different rigid normal forms due to its abundance of symmetries. The representations of the Heisenberg sphere in rigid normal form are called rigid spheres. 

We reinterpret rigid spheres as Heisenberg spheres with a Sasakian structure. In addition to the CR structure this includes the choice of an infinitesimal automorphism $Z$, i.e. a Reeb vector field the flow of which preserves the CR structure. There is an 8-parametric family of such infinitesimal automorphisms
\begin{equation}\label{ia}
Z=(p+cz+aw+2\I\bar{a}z^2+rzw )\frac{\partial}{\partial z}+ (q+2\I\bar{p}z +2(\Re c) w+2\I\bar{a}zw+rw^2 )\frac{\partial}{\partial w},
\end{equation}
where $p,a,c\in\mathbb C$, $q,r\in\mathbb R$ and $q\neq 0$.
However, different infinitesimal automorphisms may induce equivalent Sasakian structures in the following sense:
Two Sasakian manifolds $(M,Z)$ and $(M',Z')$ are equivalent if there exists a CR diffeomorphism $\Phi\colon M\to M'$ such that $\Phi_*Z=Z'$. We call two Sasakian manifolds $(M,Z)$ and $(M',Z')$ homothetic if there exists a CR diffeomorphism $\Phi\colon M\to M'$ such that $\Phi_*Z=\kappa Z'$ for some non-zero real constant $\kappa$.

The classification of rigid spheres is equivalent to the classification of homothetic Sasakian structures on the Heisenberg sphere.

Recall that Stanton's normal form of a real-analytic rigid hypersurface
$$v=\psi(z,\bar{z})$$
with $\psi(0,0)=0$, $d\psi(0,0)=0$ can be achieved by a coordinate change
\begin{align*}
z'&=\frac{\partial \psi}{\partial \bar{z}}(z,0)\\
w'&=w +\psi(z,0).
\end{align*} 

We say that two rigid hypersurfaces $M$ and $M'$ are equivalent if they are CR equivalent by a mapping $\Phi$ that takes the rigid symmetry $\frac{\partial}{\partial u}$ to the rigid symmetry $\kappa\frac{\partial}{\partial u'}$, that is $(M,\frac{\partial}{\partial u})$ and $(M', \frac{\partial}{\partial u'})$ are homothetic Sasakian manifolds.  

The Proposition below shows that the rigid normal form is an efficient tool for deciding whether two rigid hypersurfaces are equivalent.

\begin{prop}
Two real-analytic Levi-nondegenerate rigid hypersurfaces in rigid normal form $M$ and $M'$ are equivalent, that is $(M,\frac{\partial}{\partial u})$ and $(M',\frac{\partial}{\partial u'})$ are homothetic Sasakian manifolds,  if and only if they are related by a mapping
\begin{align*}
z&=cz'\\
w&= |c|^2w'.
\end{align*}
The coefficients of their normal forms are related by $\gamma'_{k\ell}=c^{k-1}\bar{c}^{\ell-1}\gamma_{k\ell}$.
\end{prop} 

{\bf Proof.} A coordinate change that pulls $\frac{\partial}{\partial u'}$ back to  $\frac{1}{\kappa}\frac{\partial}{\partial u}$  must be of the form
\begin{align*}
z&=f(z')\\
w&= \frac{1}{\kappa} w' +g(z')
\end{align*}
where $\kappa$ is a real non-zero constant and $f,g$ are functions of $z$. The absence of harmonic terms (holomorphic or antiholomorphic with respect to $z$ or $z'$) in the defining equations of $M$ and $M'$ implies that $g\equiv 0$. The absence of terms of type $(k,1)$ for $k>1$ (and their conjugates) in the defining equations implies that $f(z)=cz$ for some complex non-zero constant $c$. It is easy to see that $\frac{1}{\kappa}=|c|^2$.   \hfill $\Box$

\section{The moduli space}

In the paper  \cite{MR3515414} Ezhov and Schmalz showed that for arbitrary real numbers $\gamma_{22},\gamma_{33}$ and an arbitrary complex number $\gamma_{23}$ there is a unique rigid sphere in rigid normal form \eqref{rnf} with those coefficients, which are controlled by the parameters
\begin{align}\label{alg}
\tau&=-\frac{\gamma_{22}}{2} \\\nonumber
a&=  -\frac{\gamma_{23}}{2}\\\nonumber
\rho&=  -\frac32 \gamma_{33} + \frac94 \gamma_{22}^2.
\end{align}

Two spherical Sasakian manifolds with the parameters $(\tau,a,\rho)$ and $(\tau',a',\rho')$ are equivalent if $\tau'=\tau$, $\rho'=\rho$ and $a'=\e^{\I \phi}a$ for some $\phi\in[0,2\pi)$. Thus the moduli space of local spherical Sasakian manifolds of dimension $3$ is $\mathbb R^3_+=\{(\tau,a,\rho)\in\mathbb R^3 \mid a\ge0\}$.

The rigid spheres with parameters $(\tau,a,\rho)$ and $(\tau',a',\rho')$ are equivalent if and only if there is $c\in \mathbb C^*$ such that 
\begin{equation}\label{act}
\tau'=|c|^2\tau, \qquad a'=c\bar{c}^2a, \qquad \rho'=|c|^4\rho.
\end{equation}  

We can use the argument of $c$ to make $a$ non-negative real. The remaining action of the modulus of $c$ has the fixed point $(\tau=0,a=0,\rho=0)$, which corresponds to the Heisenberg sphere. The other orbits can be represented by points on the closed surface
\begin{equation}\label{defM}
\mathfrak M\colon \tau^4+a^{8/3}+\rho^2=1, \qquad a\ge 0,
\end{equation}
which we call the moduli space of rigid spheres. Notice that the choice of the exponents in the equation \eqref{defM} makes the expression on the left hand side homogeneous of order 8 with respect to the action \eqref{act} with positive $c$. 

This surface $\mathfrak M$ is topologically equivalent to a closed disk and it is convenient to represent it by its projection on the $\tau,\rho$-plane, that is, the closed region
 $$\tilde{\mathfrak M}\colon \tau^4+\rho^2\le1.$$

\section{The parameter space}

N. Stanton \cite{Sta90} found a class of equations of rigid spheres, 
\begin{multline}\label{sf}
\frac{1}{2r}\sin 2rv \left(1-\frac{2|b|^2\theta}{|c|^2} \right)=\frac{|z|^2 \e^{-2\theta v}}{1+4|b|^2|z|^2+ 2\I(b\bar{z}-\bar{b} z)} + \frac{|b|^2}{|c|^2}(\e^{-2\theta v}-\cos 2rv)+ \\+ \frac{\bar{b} z}{\bar{c}(1-2\I\bar{b}z)}(\e^{-2\theta v}-\e^{2\I rv})+ \frac{b\bar{z}}{c(1+2\I b\bar{z})}(\e^{-2\theta v}-\e^{-2\I rv}),
\end{multline}
which also depends on two real and one complex parameter $\theta,r,b$. The two real parameters $\theta,r$ can be combined into a complex parameter $c=r+\I \theta$. Stanton's parameters are related to our parameters by the equations
 \begin{align}\label{alg1}
\tau&=\theta -3|b|^2\\\nonumber
a&= -b (r-\I\theta + 2\I |b|^2)\\\nonumber
\rho&= -3 |b|^4 - r^2 + 2 |b|^2 \theta.
\end{align}
This raised the question whether all parameters $(\tau,a,\rho)$ can be realised by Stanton's rigid spheres and whether different parameters $(\theta,r,b)$ can give equivalent rigid spheres. In fact, the system of algebraic equations \eqref{alg1} may be unsolvable with respect to $(\theta,r,b)$ for certain triples $(\tau,a,\rho)$. 

In the paper \cite{MR3515414} Ezhov and Schmalz introduced the auxiliary parameter $\phi$, which for Stanton's rigid spheres equals $|b|^2$ but will be allowed to take negative values. Also $r^2$ will be allowed to take negative values, which renders $r$ purely imaginary.  This led to equation \eqref{sef}, which comprises all rigid spheres in terms of the parameters $(\theta,a,r,\phi)$, where $\phi$ is a real solution of the cubic equation \eqref{algseta}. Notice, that this equation has real coefficients and therefore always has at least one real solution. We can now complete the proof of Theorem \ref{th1}. Any class of equivalent rigid spheres with parameter $a\neq0$ has a unique representation by a parameter triple $(\tau',1,\rho')$, which is equivalent to the triple $$\left(\frac{\tau'}{(1+{\tau'}^4+{\rho'}^2)^4},\frac{1}{(1+{\tau'}^4+ {\rho'}^2)^{3/8}},\frac{\rho'}{(1+{\tau'}^4+{\rho'}^2)^2}\right)=(\tau,(1-{\tau}^4- {\rho}^2)^{3/8},\rho)$$
with $\tau^4+\rho^2<1$. Any class of equivalent rigid spheres with parameter $a=0$, but ${\tau'}^4+{\rho'}^2\neq 0$  has a unique representation
$$\left(\frac{\tau'}{({\tau'}^4+{\rho'}^2)^4},0,\frac{\rho'}{({\tau'}^4+{\rho'}^2)^2}\right)=(\tau,(1-{\tau}^4- {\rho}^2)^{3/8},\rho)$$
with $\tau^4+\rho^2=1$. The equivalence class with parameters $(0,0,0) $ is represented by the Heisenberg sphere.
The parameters $\theta,r,\phi$ can now be found by the formulae in Theorem \ref{th1}.

\section{The discriminant set of equation \eqref{algset}}

In this section we study the discriminant set of the cubic equation \eqref{algset}, which determines parameter $\phi$. This discriminant set is given by
$$27 a^4-18 a^2 \rho  \tau +2 a^2 \tau ^3-\rho ^3+2 \rho ^2 \tau ^2-\rho  \tau ^4=0.$$
Its projection on $\tilde{\mathfrak M}$ is illustrated in the left picture of Figure \ref{fig1} below. It is a curve that starts at $(-1,0)$ and forms a cusp as shown in the enlargement in the right picture of Figure \ref{fig1} below. It strongly osculates the boundary of $\tilde{\mathfrak M}$ at the point  $(\frac{1}{\sqrt[4]{2}},\frac{1}{\sqrt{2}})$ and terminates at the point $(1,0)$.

\begin{figure}[ht]
\includegraphics[scale=0.45]{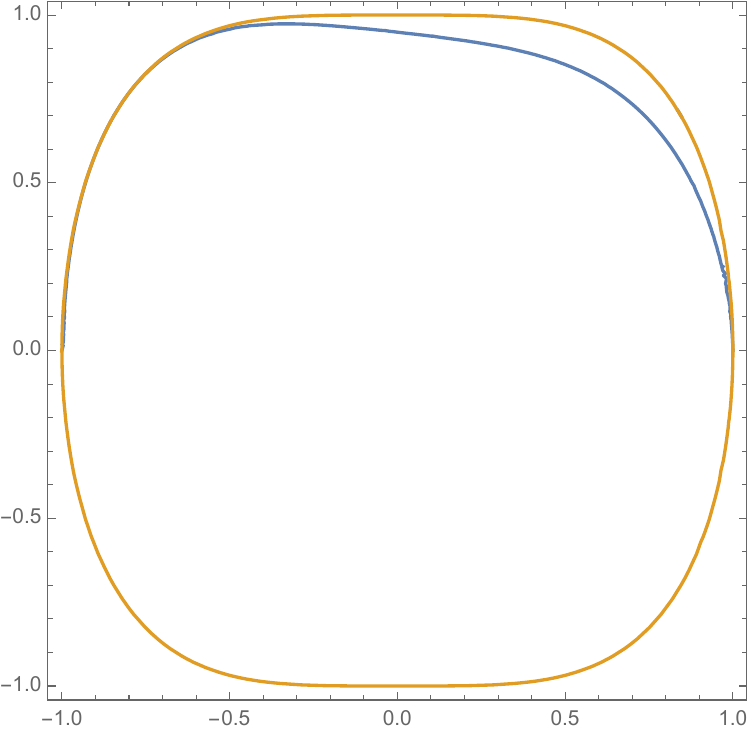}  \hfill  \includegraphics[scale=0.45]{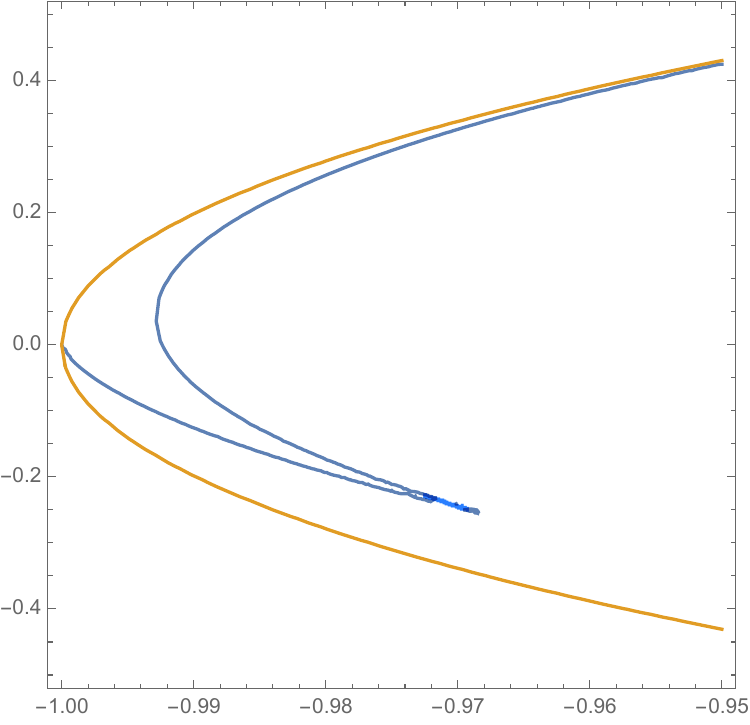} \caption{Discriminant curve}\label{fig1}
\end{figure}
\bigskip

We show now that the discriminant curve in $\tilde{\mathfrak M}$ separates the set of Stanton's rigid spheres from the remaining ones, and we will obtain a parametric representation of the curve. Stanton's rigid spheres occur for $r^2\ge 0$ and $\phi\ge 0$ (for at least one solution). Notice that  $|a|^2>0$ implies that one solution $\phi$ must be positive while the other two solutions might be both negative, complex conjugate or both positive. 

The boundary equations are $\phi=0$ or $r=0$. Consider the first case $\phi=0$. This implies $a=0$ and $r^2=-\rho \ge 0$, which gives the lower part of the boundary of  $\tilde{\mathfrak M}$.  Assume now $\phi>0$ and  $r^2=0$. Then
\begin{align*}
\rho&=(2\tau+3\phi)\phi\\
|a|^2&=(4\phi^2+4\tau\phi+\tau^2-\rho)\phi
\end{align*}
and hence
$$|a|^2=(\phi^2+2\tau\phi+\tau^2)\phi=(\phi+\tau)^2\phi.$$
That is
$$|a|=\pm(\phi+\tau)\sqrt{\phi},$$
which yields the parametric representation of the curve
\begin{equation}\label{par}
\tau=\pm\frac{|a|}{\sqrt{\phi}}-\phi \qquad \rho=\pm2|a|\sqrt{\phi}+\phi^2
\end{equation}
where $\phi\in (0,\infty)$.

In order to show that this curve coincides with the locus of real solutions of the discriminant equation we restrict to $a\neq 0$ and switch to another representation of the equivalence classes of rigid spheres by rescaling to $(\tau,2,\rho)$. For $a=2$ the discriminant equation becomes   
\begin{equation}\label{disc2}
\rho ^3-2 \rho ^2 \tau ^2+\rho  \tau ^4+72 \rho  \tau -8 \tau ^3-432=0.
\end{equation}
This choice gives a cusp with coordinates $(-3,-3)$.

Substituting \eqref{par} with $a=2$ into \eqref{disc2} confirms that the discriminant vanishes on the boundary curve of the set of Stanton's rigid spheres.  Since the discriminant locus has three real points for each $\tau<-3$,  two real points for $\tau=-3$ and one real point for $\tau>-3$ the two curves coincide. This is illustrated in Figure \ref{fig2} below.
\begin{center}
\begin{figure}[ht]
\includegraphics[scale=0.4]{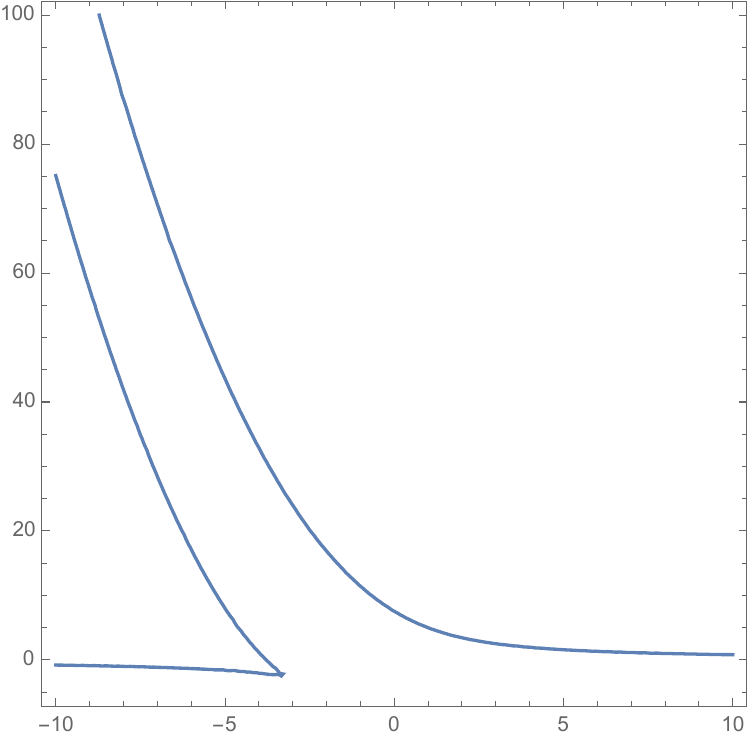} \caption{}\label{fig2}
\end{figure}
\end{center}

\begin{theo}
1. Stanton's rigid spheres are the rigid spheres represented by $(\tau,\rho)\in\tilde{\mathfrak M}$ that lie between the discriminant curve and the lower part of the boundary of $\tilde{\mathfrak M}$ (see Figure \ref{fig1}).\\
2. These are exactly the rigid spheres that possess a unique set of parameters $(\theta,a,r,\phi)$. \\
3. Two Stanton's rigid spheres are equivalent if and only if the parameters  $(\theta,b,r)$ and $(\theta', b',r')$ are related by
\begin{align*}
\theta'&=|c|^2\theta\\
b'&=\bar{c} b\\
r'&=|c|^2r.
\end{align*} 
for some non-zero complex number $c$. 
\end{theo}

{\bf Proof.} We have already proved the first two statements. Statement 3 follows from relations \eqref{alg1}, the uniqueness of the parameters $(\tau,a,\rho)$ up to the action of complex scaling and the uniqueness of the real solution $\phi$. \hfill $\Box$ \medskip

\begin{prop}\label{branch}
The cusp of the discriminant curve is a branching point for the cubic equation \eqref{algset} and there is no continuous branch of real solutions.
\end{prop}

{\bf Proof.} In the parametric equation \eqref{par} for $a=2$ the cusp corresponds to $\phi=1$. For $\phi<1$ the limit of the complex solutions on the discriminant curve is greater than the real solution and for $\phi>1$ it is smaller. If there was a continuous branch the solutions that come from the complex solutions should stay either below or above the real solution as long as we avoid the discriminant set. Since we can connect the two pieces of the discriminant curve avoiding other points of the discriminant set we conclude that no continuous branch exists. \hfill $\Box$

\section{Automorphisms of spherical Sasakian manifolds}
In this section we prove Theorem \ref{th2}.

The Heisenberg sphere with Reeb vector field (real part of) $Z= \frac{\partial}{\partial w}$ is homogeneous as a Sasakian manifold. The infinitesimal CR automorphisms
$$(p+ \I x z) \frac{\partial}{\partial z} + (q+2\I \bar{p} z) \frac{\partial}{\partial w}$$
with $p\in\mathbb C$ and $q,x\in\mathbb R$ commute with $Z$ and generate the 4-dimensional transitive affine group of Sasakian automorphisms
$$(z,u)\mapsto (\e^{\I x}(z+p), u+q-2\Im z\bar{p}).$$

Other choices of Reeb vector fields reduce to 
$$Z=(\I\tau z +aw +2\I \bar{a}z^2 +\rho zw )\frac{\partial}{\partial z} +(1+2\I \bar{a}zw +\rho w^2 )  \frac{\partial}{\partial w}.$$

The infinitesimal CR automorphisms of the Heisenberg sphere
$$\xi=(p+ xz  +\alpha w +2\I \bar{\alpha}z^2 +r zw )\frac{\partial}{\partial z} +(q+2\I \bar{p}z +2\Re x w+2\I \bar{\alpha}zw +r w^2 )  \frac{\partial}{\partial w}$$
with $p,x,\alpha\in \mathbb C$ and $q,r\in\mathbb R$ become infinitesimal Sasakian automorphisms if $[\xi,Z]=0$. This implies
\begin{align*}
\Re x&=0\\
\alpha&= a q+\I p \tau\\
r&= -\I a \bar{p}+\I \bar{a} p+q \rho\\
a \bar{p}+ \bar{a} p&=0.
\end{align*}

It follows $a=0$ or $\arg p=\arg a+\pi$.   Consider the first case. Then either $\rho -\tau ^2=0$ or $p=0$. Hence, homogeneity occurs only if $a=0$ and $\rho=\tau^2$, which can be represented by 
$$2\sinh \frac{v}{2}= e^{-\frac{v}{2}}|z|^2$$
or
$$2\sinh \frac{v}{2}= e^{\frac{v}{2}}|z|^2.$$

The first option is equivalent to
\begin{equation}\label{kae}
v=\log(1+|z|^2).
\end{equation}
Notice that the left hand side of equation  \eqref{kae} is the K\"ahler potential of the Fubini-Study metric in $\mathbb{CP}^1$ and its isometries together with the Reeb vector field $Z$ give us the automorphisms of the Sasakian manifold  \eqref{kae}, namely
\begin{align*}
z'&= \e^{\I\phi} \frac{z-\zeta}{\bar{\zeta}z+1}\\
w'&=w+q-2\I \log(1+\bar{\zeta} z) +\I \log (1+|\zeta|^2),
\end{align*}
where $q\in\mathbb R$, $\phi\in[0,2\pi)$ and $\zeta\in\mathbb C$. Globally, this Sasakian manifold can be realised as the sphere $|z|^2+|w|^2=1$ in $\mathbb C^2$, endowed with the induced metric and CR structure. The  unitary transformations in $\mathbb C^2$ induce the Sasakian automorphisms.

The second option is equivalent to
\begin{equation}\label{lob}
v=-\log(1-|z|^2)
\end{equation}
with $|z|<1$.
Here the left hand side of equation \eqref{lob} is the K\"ahler potential of the Poincar\'e-Lobachevsky metric on the unit disc and its isometries give the Sasakian automorphisms of the Sasakian manifold:
\begin{align*}
z'&= \e^{\I\phi} \frac{z+\zeta}{\bar{\zeta}z+1}\\
w'&=w+q+2\I \log(1+\bar{\zeta}
z) -\I \log (1-|\zeta|^2),
\end{align*}
where $q\in\mathbb R$, $\phi\in[0,2\pi)$ and $|\zeta|<1$. A global realisation of this Sasakian manifold is the hyperboloid $|z|^2-|w|^2=1$ in $\mathbb C^2$, endowed with the induced metric and CR structure. The  pseudounitary transformations in $\mathbb C^2$ induce the Sasakian automorphisms.

Assume now $a=0$, and  $\rho -\tau ^2\neq0$, hence $p=0$. Then the infinitesimal automorphisms  form a two-dimensional Lie algebra,  generated by $Z$ and the rotation $\I z \frac{\partial}{\partial z}$.

Finally, assume $a\neq 0$. Set $p=\I  t a$, with $t\in\mathbb R$. Then 
\begin{align*}
x&=\I (t(\rho-\tau^2)+q\tau)\\
\alpha&=a(q-t\tau)\\
r&=-2t |a|^2 +q\rho.
\end{align*}
This results in a 2-dimensional Lie algebra spanned by   $Z$ and
$$(\I a+  \I(\rho-\tau^2)z  -a \tau w -2\I \bar{a}\tau z^2 -2|a|^2 zw )\frac{\partial}{\partial z} +(2\bar{a}z- 2\I \bar{a} \tau zw -2|a^2| w^2 )  \frac{\partial}{\partial w}.$$

\section{Examples}
{\bf Example 1.} First we consider the case $a=0$. Then the equation on $\phi$ takes the simple form
$$\phi^3+ \tau\phi^2+ \frac{\tau^2-\rho}{4}\phi=0.$$
It has always the solution $\phi=0$ but also the solutions
$$\phi=-\frac{\tau}{2}\pm \frac{\sqrt{\rho}}{2}$$
which are real if $\rho\ge 0$. We get the following options:

\begin{align}
\phi=&0,& \I r=&\sqrt{\rho},& \theta=&\tau\\
\phi=&\frac12(-\tau+\sqrt{\rho}),& \I r=&\frac12(\tau+\sqrt{\rho}),& \theta=&\frac12(-\tau+3\sqrt{\rho})\\
\phi=&\frac12(-\tau-\sqrt{\rho}),& \I r=&\frac12(\tau-\sqrt{\rho}),& \theta=&\frac12(-\tau-3\sqrt{\rho})
\end{align}
In each case the resulting rigid sphere is
$$\frac{\sinh 2\sqrt{\rho}v}{ 2\sqrt{\rho}}=\e^{-2\tau v}|z|^2,$$
where each equivalence class can be represented by the pair $(\tau, \rho=\sqrt{1-\tau^4})$ for $\tau\in[-1,1]$. The homogeneous manifolds $v=\pm\log(1\pm|z|^2)$ from the previous section belongs to this family with $\tau=\pm\frac12$ and $\rho=\frac14$. 
\medskip

{\bf Example 2.} The rigid sphere with $\tau=\rho=0$ and $a=2$ is in the region of Stanton's rigid spheres. 
In this case $\theta=3\phi$, $r^2=3\phi^2$. From
\begin{align*}
2=&a=-b(\sqrt{3}-\I)\phi\\
4=&|a|^2=4\phi^3
\end{align*}
it follows $\phi=1$, $\theta=r=3$, $b=\frac{-2}{\sqrt{3}-\I}$.\medskip

{\bf Example 3.} The cusp occurs for $a=2$, $\tau=\rho=-3$, which gives the triple solution $\phi=1$ and hence $r=0$, $\theta=4$, $b=-\I$. This example appears in Stanton's list in \cite{Sta90} as (6.17).


\begin{bibdiv}
\begin{biblist}
\bib{BER}{book}{
   author={Baouendi, M. Salah},
   author={Ebenfelt, Peter},
   author={Rothschild, Linda Preiss},
   title={Real submanifolds in complex space and their mappings},
   series={Princeton Mathematical Series},
   volume={47},
   publisher={Princeton University Press},
   place={Princeton, NJ},
   date={1999},
   pages={xii+404},
   isbn={0-691-00498-6},
   review={\MR{1668103 (2000b:32066)}},
}
\bib{MR3515414}{article}{
   author={Ezhov, Vladimir},
   author={Schmalz, Gerd},
   title={Explicit description of spherical rigid hypersurfaces in
   $\mathbb{C}^2$},
   journal={Complex Anal. Synerg.},
   volume={1},
   date={2015},
   number={1},
   pages={Paper No. 2, 10},
   issn={2197-120X},
   review={\MR{3515414}},
   doi={10.1186/2197-120X-1-2},
}
\bib{MR4021084}{article}{
   author={Isaev, Alexander},
   author={Merker, Jo\"{e}l},
   title={On the real-analyticity of rigid spherical hypersurfaces in
   $\mathbb{C}^2$},
   journal={Proc. Amer. Math. Soc.},
   volume={147},
   date={2019},
   number={12},
   pages={5251--5256},
   issn={0002-9939},
   review={\MR{4021084}},
   doi={10.1090/proc/14724},
}
\bib{Sta90}{article}{
   author={Stanton, Nancy K.},
   title={A normal form for rigid hypersurfaces in ${\bf C}^2$},
   journal={Amer. J. Math.},
   volume={113},
   date={1991},
   number={5},
   pages={877--910},
   issn={0002-9327},
   review={\MR{1129296}},
   doi={10.2307/2374789},
}
\end{biblist}
\end{bibdiv}
\end{document}